\documentclass{amsart}

\usepackage{amsmath,amssymb,amsthm,amsfonts, mathrsfs}
\usepackage[all]{xy}
\usepackage{enumerate}
\usepackage{graphicx}

\newtheorem{theorem}{Theorem}[section]

\newtheorem{proposition}[theorem]{Proposition}

\newtheorem{definition}[theorem]{Definition}

\begin{document}

\title{Extension properties of orbit spaces of proper  actions  revisited}

\author{Sergey Antonyan\\[3pt]}

\address{Departamento de  Matem\'aticas,
Facultad de Ciencias, Universidad Nacional Aut\'onoma de M\'exico,
 04510 M\'exico City, M\'exico.}
\email{antonyan@unam.mx}

\begin{abstract} Let  $G$ be a locally compact Hausdorff group. We   study orbit spaces   of equivariant absolute neighborhood extensors ($G$-${\rm ANE}$'s) in the class
 of all proper $G$-spaces  that are metrizable by a $G$-invariant metric. 
 We  prove that if a proper $G$-space $X$ is a  $G$-${\rm ANE}$ and all  $G $-orbits in $X$ are metrizable, then  the $G$-orbit space  $X/G$ is an {\rm ANE}. If $G$ is a Lie group and 
 $H$  is a closed  normal subgroup of $G$, then  the $H$-orbit space  $X/H$ is a $G/H$-{\rm ANE}. 
\end{abstract}


\maketitle
\markboth{SERGEY ANTONYAN}{ORBIT SPACES  OF $G$-${\rm ANE}$'S}

\section {Introduction}

The  main  purpose  of this note is to prove the following two theorems.

\begin{theorem}[Orbit space theorem]\label{T:0}
Let $G$ be a locally compact Hausdorff group and $X$   a proper $G$-space such that  all  $G$-orbits in $X$ are metrizable. 
 If  $X$ is a $G$-${\rm ANE}$, then the $G$-orbit space  $X/G$ is an {\rm ANE}. 

\end{theorem}

\smallskip

This theorem was first  proved in  \cite[Theorems 3.11]{ant:99} for   an almost connected acting group $G$ (i.e., the space of connected components  of $G$ is compact) and a phase space $X$ with  a paracompact orbit space $X/G$.  In \cite[Theorem 6.4]{ant:jap}      a proof  of this theorem was provided without any additional restrictions. In that proof         
      the following affirmation that we state here in the form of a proposition was used.

 \begin{proposition}\label{retractdirect}
 Let $G$ be a topological group and $K$ a  compact subgroup of  $G$. If $S$ is a $K$-space, then 
 $({G \times_K S})/G$ is homeomorphic to a retract of  $({G \times_K S})/K$.
 \end{proposition}
 
But the argument for the proof of this statement  given in \cite[Theorem 6.4]{ant:99},  unfortunately,  works correctly only for an abelian acting group $G$.  Namely,  in  that  proof 
 the formula  $(G \times_K S)/G\cong G/K\times S/K$ was used,  which however is true only for an abelian group $G$ (see \cite[Proposition 2]{ant:99}).
 
  Below, in Section \ref{Prop}, we provide a simple proof of this proposition for any topological  group $G$, thus filling the gap in the proof  of \cite[Theorem 6.4]{ant:99}.
 
 \medskip
  
 The second theorem is the following.
 
   \begin{theorem}[The case of Lie group actions]\label{T:00} Let $G$ be a Lie  group,  $H$  a closed normal subgroup of $G$, and $X$   a proper $G$-space. If 
 $X$ is a $G$-${\rm ANE}$, then the $H$-orbit space  $X/H$ is a $G/H$-${\rm ANE}$. 

\end{theorem} 

 In \cite[Theorem 1.1]{ant:proper}      a proof  of this theorem was given even for  any locally compact acting group $G$. Again, in that proof a formula was used (see \cite[formula (3.3)]{ant:proper} which is correct only for abelian groups. Below we will provide a very short proof of this theorem in the case of Lie group actions, which is practically the most important case.
 This proof is based on the following our result  proved 
 in \cite[Proposition 4.1]{ant:proper}.
  
  \begin{proposition}\label{P:1} Let $G$ be a Lie  group, $K$ a compact   subgroup of $G$, and $S$  a $K$-space. Then $S$ is a neighborhood $K$-equivariant retract of the twisted product $G\times_K S$.
\end{proposition}    
        \medskip

Recall that the orbit space problem was  posed in \cite[Question 4]{ant:aspects}. It has been  solved first 
 in \cite[Theorem 8]{ant:88} where it was established that if $G$ is a compact metrizable group and $X$  a $G$-A(N)R, then the orbit space $X/G$ is an A(N)R. 
 This result was widely applied in the study of the topology of  Banach-Mazur compacta  (see \cite{ant:98}, \cite{ant:00}, \cite{ant:03}).  Other applications  can be found in \cite{ant:88}, \cite{bas:94} and \cite{zhur:90}.

 \medskip
Before passing to the details of the proofs it is convenient to recall some auxiliary notions and results.
 
\medskip

\section {Some basic definitions and auxiliary results}

Throughout the paper the letter $G$ will denote a locally compact Hausdorff group unless otherwise is stated; by $e$ \ we  denote the unity of $G$. 

All topological spaces  are assumed to be Tychonoff (= completely regular and Hausdorff). 
The basic ideas and facts of the theory of $G$-spaces or topological
transformation groups can be found in Bredon \cite{bre:72} and in  Palais \cite{pal:60}.  
Our basic references on  proper group actions are   Palais~\cite{pal:61} and Abels   \cite{abe:78}.        
For the  equivariant theory of retracts the reader can see, for instance, \cite{ant:87},  \cite{ant:88} \cite{ant:99}, \cite{ant:jap} and \cite{ant:proper}. 

For the  convenience of the reader we recall, however,  some more  special definitions and facts.

\smallskip

Here we deal with $G$-spaces. If $X$ and $Y$ are two $G$-spaces then a continuous map $f:X\to Y$ is called a $G$-map, if $f(gx)=gf(x)$ for all $x\in X$ and $g\in G$. If a $G$-map is a homemorphism then it is called a $G$-homeomorphism.

If $X$ is a $G$-space and $H$ a subgroup of $G$ then, for a subset $S\subset X$, \ $H(S)$ denotes the $H$-saturation of $S$, i.e., $H(S)$= $\{hs | \ h\in H,\ s\in S\}$. In particular, $H(x)$  denotes the $H$-orbit $\{hx\in X | \ h\in H \}$ of $x$.  The quotient space  of all $H$-orbits is called the $H$-orbit space and denoted by $X/H$. 

If $H(S)$=$S$, then $S$ is said to be an $H$-invariant set. A  $G$-invariant set will simply be called an invariant set.  

For a closed subgroup $H \subset G$, by $G/H$ we will denote the $G$-space 
of cosets $\{gH | \ g\in G\}$ under the action induced by left translations.

If $X$ is a $G$-space and $H$  a closed normal subgroup of $G$, then the $H$-orbit space $X/H$  will always be regarded as a  $G/H$-space endowed with the following action of the group $G/H$:
$(gH)*H(x)=H(gx), \ \ \text{where} \ \ gH\in G/H, \ H(x)\in X/H$.

For any $x\in X$, the subgroup   $G_x =\{g\in G   \mid  gx=x\}$ is called  the stabilizer (or stationary subgroup) at $x$.

\smallskip

 Let  $X$ be  a $G$-space. Two subsets  $U$ and $V$ in  $X$  are called  thin relative to each other   \cite[Definition 1.1.1]{pal:61},  if the set 
$\langle U,V\rangle=\{g\in G | \ gU\cap V\ne \emptyset\}$
    has  a compact closure in $G$. 
   A subset $U$ of a $G$-space $X$ is called {\it small}, \ if  every point in $X$ has a neighborhood thin relative to $U$. A $G$-space $X$ 
is called  {\it  proper} (in the sense of R. Palais),  if   every point in  $X$ has a small neighborhood. We refer to the seminal paper of R. Palais \cite{pal:61} for further information about proper $G$-spaces.

\smallskip

In the present  paper we are especially interested in the class $G$-$\mathcal M$ of  all  metrizable proper $G$-spaces  that admit a compatible $G$-invariant metric. 
It is well-known that, for $G$ a compact group, the class $G$-$\mathcal M$
coincides with the class of {\it all} metrizable $G$-spaces (see \cite[Proposition 1.1.12]{pal:60}).
A  fundamental result of  R.~Palais \cite[Theorem~4.3.4]{pal:61}
states that if $G$ is a Lie group, then  $G$-$\mathcal M$ includes all {\it separable}, metrizable  proper $G$-spaces.

 Let us recall the definition of a twisted product $G/H\times_K S$, where $H$ is a closed normal subgroup of $G$,  $K$  any closed subgroup of $G$,  and   $S$  a $K$-space. 
 
$G/H\times_KS$ is the orbit space of the $K$-space $G/H\times S$, where $K$ acts on the Cartesian product $G/H\times S$ by $k(gH, s)=(gk^{-1}H, ks)$. Furthermore, there is a natural action of $G$ on $G/H\times_K S$ given by $g^\prime[gH, s]=[g^\prime gH, s]$, where $g'\in G$ and $[gH, s]$ denotes the  $K$-orbit of the point $(gH, s)$ in $G/H\times S$. The twisted products of the form $G\times_KS$ (i.e., when $H$ is  the trivial subgroup of $G$) are  of a particular interest in the theory of transformation groups (see \cite[Ch.~II, \S~2]{bre:72}). 

 \smallskip

A $G$-space  $Y$ is called an equivariant absolute neighborhood extensor  for the class  $G$-$\mathcal M$ 
(notation: $Y\in G$-{\rm ANE})  if,   for any  $X\in G$-$\mathcal M$  and any closed invariant subset $A\subset X$, every    $G$-map $f:A\to Y$ admits a $G$-map $\psi\colon U\to Y$ defined on 
  an invariant neighborhood $U$ of $A$ in $X$ such that $\psi|_A= f$. If, in addition, one  can always take $U=X$, then we say that $Y$ is an equivariant absolute extensor  for  $G$-$\mathcal M$ 
  (notation: $Y\in G$-AE). The map $\psi$ is called a $G$-extension of $f$.

\smallskip

Let us recall  the well known definition of a slice \cite[p.~305]{pal:61}:

\begin{definition}\label{D:21} Let $X$ be a $G$-space and   $H$   a closed  subgroup of $G$. An $H$-invariant  subset $S\subset X$ is called an  $H$-slice in $X$, if $G(S)$ is open in $X$ and there exists  a $G$-map $f:G(S)\to G/H$ such that $S$=$f^{-1}(eH)$.  The saturation $G(S)$ is  called   a {\it tubular} set and $H$ is called a slicing group. 

  If  $G(S)=X$, then we say that $S$ is {\it a global} $H$-slice for $X$.   
\end{definition}
\medskip

The following  result of R. Palais \cite[Proposition 2.3.1]{pal:61} plays  a central role in the theory of topological transformation groups.

\begin{theorem}[Slice Theorem]\label{T:Sl} Let $G$ be a Lie group, $X$ be a proper $G$-space and $x\in X$. Then there exists a $G_x$-slice $S\subset X$ such that $x\in S$.
\end{theorem}

In our proofs we will also need  the following approximate version of the Slice Theorem  proved in  \cite[Theorem 3.6]{ant:jap} (see also \cite[Theorem 6.1]{ant:dikr})  which is valid for any locally compact group.

\begin{theorem}[Approximate Slice Theorem] \label{T:331} Let $G$ be any group, $X$  a proper $G$-space and  $x\in X$. Then for any  neighborhood $O$ of $x$ in $X$,   there exist  a   compact  large subgroup $K$ of $G$ with $G_x\subset K$, and a $K$-slice  $S$ such that  $x\in S\subset O$. 
\end{theorem}

Recall that here a subgroup $K\subset G$ is called  {\it large}, if there exists a closed  normal subgroup $N\subset G$ such that
$N\subset K$ and $G/N$ is a Lie group.

In the context of equivariant extension properties the  notion of a  large subgroup was first singled out in \cite{ant:94} (for compact groups)  and in \cite{ant:99} (for locally compact groups). Although  some geometric characterizations of this notion were available much earlier (see \cite[Section 3]{ant:dikr} and  the literature cited there), new characterizations were given in \cite[Proposition 6]{ant:99}, \cite[Proposition 3.2]{ant:jap} and \cite[Theorem 5.3]{ant:dikr} through equivariant extension properties.

\medskip

One of the strong properties of large subgroups is expressed in the following

\begin{proposition}[{\cite[Proposition 3.4]{ant:jap}}]\label{P:large} Let $K$  be a compact  large  subgroup of $G$, and $X$  a  $G$-${\rm ANE}$ (respectively, a $G$-$AE$). Then  $X$ is  a $K$-${\rm ANE}$ (respectively, a $K$-$AE$).
\end{proposition}

\medskip
  
  The following proposition is well known  (see, e.g. \cite[Lemma 3.5]{abe:78}).
  
  \begin{proposition}\label{twist} Let  $H$ be a compact subgroup of $G$, $X$  a proper $G$-space and  $S$ a global $H$-slice of $X$.  Then  the  map  $\xi:G\times_H S\to X$ defined by $\xi([g, s])=gs$ is a $G$-homeomorphism.
\end{proposition}

The following two results are also used in our proofs.

\medskip

\begin{theorem}[\cite{ant:99}]\label{T:21}
Let $G$ be a compact group and  $H$  a closed normal subgroup of $G$. Suppose $X$ is a $G$-space such that all  $H$-orbits in $X$ are metrizable. 
 If  $X$ is a $G$-${\rm ANE}$ (respectively, a $G$-$AE$),  then the $H$-orbit space $X/H$ is a $G/H$-${\rm ANE}$ (respectively, a $G/H$-$AE$). 

\end{theorem}

We refer to  \cite[Theorem 1]{ant:99} for the  details.

\medskip

The following  equivariant version of   Hanner's open union theorem \cite[Theorem 19.2]{ha:52} is proved in \cite[Corollary 5.7]{ant:jap}. A  short and beautiful proof of Hanner's theorem 
was given by J. Dydak \cite[Corollary 1.5]{dy:pp}.

\begin{theorem}[\cite{ant:jap}]\label{T:Un} Let $Z\in G$-$\mathcal M$.
 If a  $G$-space $Y$ is the  union of a  family of  invariant open 
$G$-${\rm ANE}(Z)$ subsets $Y_\mu\subset Y$, $\mu\in \mathcal M$, then $Y$ is a $G$-${\rm ANE}(Z)$. 
\end{theorem}

\smallskip

\section{Proof of Proposition \ref{retractdirect}}\label{Prop} 

 For every $[g, s]\in G \times_K S$ we will denote by $[g, s]_G$ the $G$-orbit in the $G$-space $G \times_K S$. Similarly, $[g, s]_K$ will denote the $K$-orbit of $[g, s]$ in $G \times_K S$.  
  
 Define the map   $\iota: (G\times_K S)/G  \to (G \times_K S)/K$ by the formula $\iota: [g, s]_G\mapsto [e, s]_K$, where $e\in G$ is the unit element. This map is well defined since for any 
 $k\in K$ one has 
 $$\iota : [gk^{-1}, ks]_G\mapsto [e, ks]_K=[k, s]_K= (k[e, s])_K= [e, s]_K.$$
 The continuity of $\iota$ is also  evident. Indeed, denote by $i$  the composition of the following continuous maps:
   $$G\times S   \to G \times S\to G\times_K S\to  (G \times_K S)/K. $$ 
    $$(g, s)\mapsto (e, s)\mapsto [e, s]\mapsto [e, s]_K.$$   
    Observe that  $i:G\times S  \to (G \times_K S)/K$ is constant on the $K$-orbits of the $K$-space $G\times S$, and hence due to compactness of $K$, it induces a continuous map $j:G\times_K S  \to (G \times_K S)/K$, $j([g, s])=    [e, s]_K$. In turn, $j$ is constant on the $G$-orbits of the $G$-space $G \times_K S$, and thence, it induces a continuous map which is exactly $\iota$

             \[ \xymatrix{G \times S \ar[r]^i \ar[d] \ar[d]_p\ar[r]  &  (G \times_K S)/K \ar@/^/[d]^ r\\
              G \times_K S \ar [r]^q\ar@{-->}[ur]^j   &  (G \times_K S)/G  \ar[u]^\iota,
             }   \]             
             where $p$ and $q$ are the orbit maps.          Thus, $\iota$ is continuous.
             
              Next we define a continuous map  $r: (G\times_K S)/K  \to (G \times_K S)/G$ by the formula $r: [g, s]_K\mapsto [g, s]_G$.

        Observe that for every $[g, s]_G\in  (G \times_K S)/G$,
    \[ r {\iota}([g, s]_G) = r([e, s]_K) = [e, s]_G =  (g[e, s])_G=[g, s]_G. \]

 Thus, $r$ is the right inverse of ${\iota}$. This implies that ${\iota}$ is an embedding and its  image ${\iota}((G \times_K S)/G)$ is a retract of     $(G \times_K X)/{K}.$ Hence, $(G \times_K S)/G$ is homemorphic to a retract of $(G \times_K X)/{K},$ as requried. 
\qed

   \
  
  \section{ Proof  of Theorem \ref{T:0}}
  
   By Theorem~ \ref{T:331}, $X$ has  an open invariant cover by tubular sets of the form $G(S)$, \  where each $S$ is a $K$-slice with  the slicing group $K$  a   compact large  subgroup of $ G$. Then the orbit space $X/G$ is the  union of its open subsets of the form $G(S)/G$. According to  
Hanner's open union theorem  \cite[Theorem 19.2]{ha:52} or \cite[Corollary 1.5]{dy:pp} (see also  
Theorem \ref{T:Un}), it suffices to show that each $G(S)/G$ is an {\rm ANE}.

To this end, we first observe that  each $G(S)$ is $G$-homeomorphic to the twisted product $G\times_{K} S$ (see Proposition \ref{twist}). This implies that $G(S)/G$ is homeomorphic to $(G\times_{K} S)/G$.
Since $X\in G$-{\rm ANE}, the tubular set  $G(S)$, being  an  open invariant subset of  $X$, is itself  a $G$-{\rm ANE}. Thus, $G\times_K S$ is a $G$-{\rm ANE}. Since the slicing group $K$  is a compact large subgroup of $G$, one can apply Proposition \ref{P:large}, according to which 
$G\times_{K} S$ is a $K$-{\rm ANE}. 
Each $K$-orbit in $X$ is contained in a $G$-orbit, and hence, is metrizable. Since $K$ is compact,  Theorem \ref{T:21} implies that
$(G\times_{K} S)/K$ is an {\rm ANE}. By Proposition \ref{retractdirect}, $(G\times_{K} S)/G$ is homeomorphic to a retract of $(G\times_{K} S)/K$, and hence, is itself  an {\rm ANE}. Consequently, $G(S)/G$ is an {\rm ANE},  as required.

\qed
  
  \

  \section{ Proof  of Theorem \ref{T:00}}

 By Theorem~ \ref{T:Sl},
 $X$ has  an open invariant cover by tubular sets of the form $G(S)$, \  where each $S$ is a $K$-slice  with the slicing group  $K$ a compact   subgroup of $G$.  Then the $G/H$-space $X/H$ is the  union of its open $G/H$-invariant subsets of the form $G(S)/H$. According to  Theorem \ref{T:Un}, it suffices to show that each $G(S)/H$ is a $G/H$-{\rm ANE}.

To this end, we first observe that  each $G(S)$ is $G$-homeomorphic to the twisted product $G\times_{K} S$ (see Proposition \ref{twist}). This yields that $G(S)/H$ is $G/H$-homeomorphic to 
$(G\times_ K S)/H$.
Next, since $X\in G$-{\rm ANE}, the tubular set  $G(S)$, being  an  open invariant subset of  $X$, is itself  a $G$-{\rm ANE}. Thus, $G\times_{K} S$ is a $G$-{\rm ANE}. Since $G$ is a Lie group we infer that $K$ is a compact large subgroup of $G$. Then one can apply Proposition \ref{P:large}, according to which 
$G\times_{K} S$ is a $K$-{\rm ANE}. 
By Proposition 1, $S$ is a $K$-equivariant retract of  $G\times_{K} S$, and hence, is  a $K$-{\rm ANE}.

Further, one has the following $G$-homeomorphism $(G\times_K S)/H\cong G/H\times_ K S$ (see \cite[Proposition 3.3]{ant:proper}). Since $S\in K$-{\rm ANE},  it then follows that the twisted product $G/H\times_ K S$ is a $G/H$-{\rm ANE} (see \cite[Proposition 3.3]{ant:proper}). This yields that  $(G\times_ K S)/H\in G/H$-{\rm ANE}, and since, $G(S)/H$ is $G/H$-homeomorphic to 
$(G\times_ K S)/H$, we conclude that $G(S)/H\in G/H$-{\rm ANE}, as required.
 \qed

\smallskip

\bibliographystyle{amsplain}

\bibliography{triquot}
\end{document}